\documentclass[cup7b]{cupbook}
\usepackage{datpatch}
\usepackage{amsmath}

\usepackage{graphics}
\DeclareFontEncoding{OT2}{}{} % to enable usage of cyrillic fonts
  \newcommand{\textcyr}[1]{%
    {\fontencoding{OT2}\fontfamily{wncyr}\fontseries{m}\fontshape{n}%
     \selectfont #1}}
\newcommand{\Sha}{{\mbox{\textcyr{Sh}}}}
\newtheorem{prop}{Proposition}[section]
\newtheorem{question}[prop]{Question}

\begin{document}

\author[Mark Watkins]{Mark Watkins \\ University of Bristol}
\chapter{Rank distribution in a family of cubic twists}

\vspace{-12pt}
\begin{abstract}
In 1987, Zagier and Kramarz published a paper in which they presented
evidence that a positive proportion of the even-signed cubic twists
of the elliptic curve $X_0(27)$ should have positive rank.
We extend their data, showing that it is more likely
that the proportion goes to zero.
\end{abstract}

\section{Introduction}
Let $E_m$ be the elliptic curve defined by the equation $x^3+y^3=m$,
which is isomorphic to $y^2=x^3-432m^2$.
The case of $m=1$ is the curve~$X_0(27)$, and the cubefree positive
$m$-values correspond to cubic twists.

These equations have a long history, dating back to Fermat.
An early study was done by Sylvester \cite{sylvester} in 1879-80,
and another voluminous study in 1951 by Selmer \cite{selmer}.
In between these two, Nagell \cite[p.14]{nagell} proved sundry
results concerning non-solvability in many cases.
In the late 1960s, Stephens \cite{stephens1,stephens2} did
numerical experiments with these curves with respect to the
then-new Birch--Swinnerton-Dyer conjecture. Zagier and Kramarz
\cite{zk} did a large numerical experiment in the~1980s,
which led them to suggest that a positive proportion of the curves
have rank 2 or greater. The best results in this regard appear
to be due to Mai \cite{mai1,mai2}, who showed that, assuming the
Parity Conjecture, for every $\epsilon>0$
at least $c_\epsilon T^{2/3-\epsilon}$ of the cubefree even twists
up to $T$ have rank 2. Elkies and Rogers \cite{er} have recently found
that the curve $$x^3+y^3=13293998056584952174157235$$ has rank at least 11.
We shall mainly be concerned with rank 2 cubic twists and in
extending the numerical data of \cite{zk}, showing that the purported
positive proportion does not seem to persist.
We also consider questions of the distribution of the size
of the Tate--Shafarevitch groups attached to these curves,
comment on effects stemming from the arithmetic of~$m$,
consider similar questions for quartic twists of $X_0(32)$,
and discuss random matrix models for these.

We briefly review how to compute the central $L$-value of $E_m$.
The first consideration is the sign of the functional equation,
which was computed by Birch and Stephens \cite{bs}.
This is defined by $\epsilon=\prod_p \epsilon_p$ where for $p\neq3$
we have that $\epsilon_p=\bigl(\frac{p}{3}\bigr)$ if~$p|m$
and $\epsilon_p=+1$ if $p$ does not divide~$m$.
For $p=3$, we have that $\epsilon_3=+1$ if $m\equiv \pm 1$ (mod~9) or $3\|m$,
and $\epsilon_3=-1$ otherwise.
Next, there is the conductor $N=\prod_p N_p$ where for $p\neq3$
we have that $N_p=p^2$ if $p|m$ and $N_p=1$ otherwise,
while for $p=3$ we have that $N_3=3^5$ if $3|m$,
that $N_3=3^2$ if $m\equiv\pm 2$ (mod~9), and $N_3=3^3$ otherwise.
There are also Tamagawa numbers and considerations for the real
period~$\Omega$; the effects of these are given in the last section
(see also Table~1 of~\cite{zk}).

When $\epsilon=+1$, the central $L$-value is given by
$$L(E_m,1)=2\sum_n \frac{a_m(n)}{n} e^{-2\pi n/\sqrt N},$$
where the conductor~$N$ is defined as above,
and the $a_m(n)$ can be computed as follows.
For primes $p\not\equiv 1$ (mod~3) and primes~$p|3m$, we define~$a_m(p)=0$.
Given a prime~$p\equiv 1$ (mod~3), the set
$$A_p=\{a|\>\> a\equiv 2\, ({\rm mod}~3),\, a^2+3b^2=4p
\>\>{\rm for\,\, some\,\,} b\in{\bf Z}\}$$ has 3 elements.
For such a prime we define $a_1(p)$ to be the unique element in $A_p$
for which $3|b$. We then define $a_m(p)$ uniquely
by the conditions $a_m(p)\equiv m^{(p-1)/3}a_1(p)$ (mod~$p$)
and $a_m(p)\in A_p$ (this second condition is equivalent to
$|a_m(p)|<2\sqrt p$ for $p>13$ and not $p\ge 13$ as \cite{zk}~claims).
Having defined $a_m(p)$ for all primes~$p$,
we extend it to prime powers via the Hecke relations,
and then to all positive integers via multiplicativity.
In order to approximate $L(E_m,1)$ well,
we need to use about $C\sqrt N$ coefficients for some constant~$C$.
When $\epsilon=-1$,
the series for $L'(E_m,1)$ has the exponential function
replaced by an exponential integral --- we did not deal with this case
(\cite{zk}~considered it for~$m\le 20000$)
since the exponential homomorphism can be computed rapidly
more readily than the exponential integral --- for the latter,
local power series would likely be useful.
Lieman \cite{lieman} has shown that the values
of $L(E_m,1)$ are the coefficients of a metaplectic form as was suggested
in~\cite[3.1]{zk}, but this does not seem useful for computational
purposes. We did not try to use the conditions given by Rodriguez Villegas
and Zagier~\cite{RVZ}, and cannot comment on their computational efficacy.

\section{Numerical data}
Applying the above method for the cubefree $m\le 10^7$ with $\epsilon=+1$,
we find that about 17.7\% of the twists have vanishing central $L$-value.
This is to be compared to 23.3\% for the $m\le 70000$,
and 20.5\% for $m\le 10^6$.
If we take the best linear fit to a log-log regression, we find that
the number of twists up to $x$ with vanishing central $L$-value
appears to grow like $x^{0.935}$. Heuristic models involving the
expected size of $\Sha$ as in \cite{zk} imply that the growth should be
more like $x^{5/6}$. Stronger models such as those in \cite{ckrs}
imply this should be more like $Bx^{5/6}(\log x)^C$
for some constants $B$~and~$C$; in the last section we make remarks
about what random matrix theory implies about~$C$.

\begin{figure}[hbtp]
\begin{center}
\scalebox{0.86}{\includegraphics{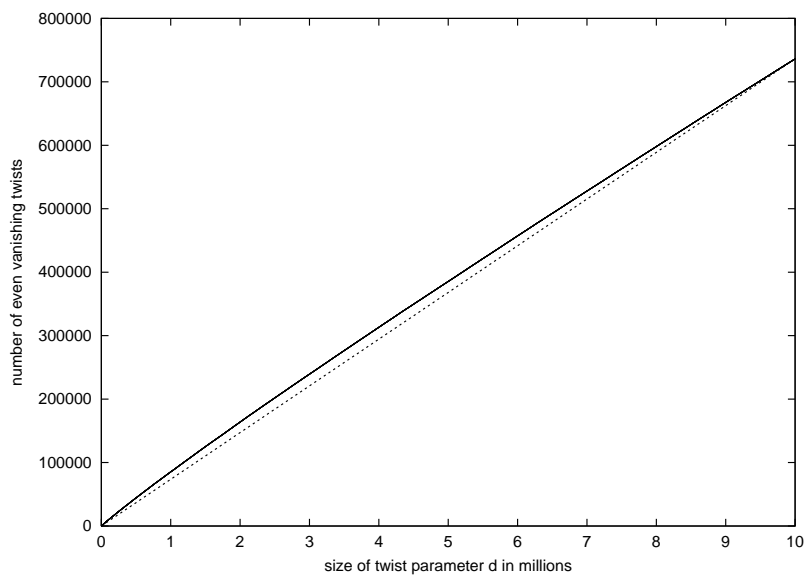}}
\caption{Number of even vanishing cubic twists of $X_0(27)$
compared to a (dotted) straight line.}
\end{center}
\end{figure}

There is also the question of arithmetic effects of~$m$.
Only 6.1\% of the prime~$m$ in the above range
have vanishing central $L$-value, while 11.3\% of the $m$
with two prime factors do, and 17.1\% of the~$m$ with three prime factors.
The number grows to 24.5\% for four prime factors, and 35.3\% for five
prime factors, and is 51.4\% for six or more prime factors.
However, each of these percentages is about 20\% lower than
the comparative value when considering only the $m\le 10^6$.
So even if we restrict to prime~$m$ we expect that the proportion
of twists with non-vanishing central $L$-value tends to zero.
Note in this context that 3-descent can tell us much about the rank
when we limit the number of prime factors of~$m$ (see \cite{cassels}).
For instance, when $m$ is prime and $E_m$ has even functional equation,
we know that $m\equiv 1,2,5$ (mod~9), and the rank is zero in
the latter two cases. Thus the 6.1\% of above might be re-interpreted
as 18.3\% of the cases where descent considerations do not force the
rank to be zero. Using the results of \cite{nagell}, we could similarly
derive such results when $m$ has two prime factors. Also, one can recall
that Elkies (see \cite{elkies}) has proven that the rank is exactly 1 for
primes $m\equiv 4,7$ (mod~9); here in fact the conjecture is that the
same is true for $m\equiv 8$ (mod~9). We return to such considerations
below when we discuss random matrix models.

We next make some comments about how often various $|\Sha|$-values occur.
Zagier and Kramarz found that 26.3\% of the even twists for $m\le 70000$
have rank 0 and trivial $\Sha$, while we find the percentage to be 18.8\%
for $m\le 10^6$ and 14.1\% for $m\le 10^7$. Indeed, already in \cite{zk}
this percentage was noted to be diminishing.
More interesting might be how often a given prime divides~$|\Sha|$,
under the restriction to rank 0 twists.
For instance, 32.4\% of the even rank 0 twists with $m\le 70000$
have 3 dividing~$|\Sha|$. This number increases to 40.1\% for $m\le 10^6$,
and is 45.3\% for $m\le 10^7$. The heuristics of \cite{delaunay} imply a
number more like 36.1\%. There is a strong arithmetic impact from~$m$,
as for prime~$m$ the percentage for $m\le 10^7$ is only 5.8\%.
However, this last datum should probably be considered anomalous
because of the special role that 3 plays in the cubic twists.

Similarly, 2~divides $|\Sha|$ about 45.7\% of the time for
even rank 0 twists with $m\le 10^7$, while only 42.1\% of the time
for $m\le 10^6$ and 35.5\% of the time for $m\le 70000$. Here Delaunay
predicts 58.1\%. Here prime~$m$ are {\it more} likely to cause 2-divisibility
of~$|\Sha|$, with the percentage here for \hbox{$m\le 10^7$} being 53.5\%.
As \cite{zk} notes, the expectation is that $|\Sha|$ should be of
size $m^{1/3}\approx N^{1/6}$ for these cubic twists, larger than
the expected $N^{1/12}$ in the general case. For 5-divisibility of~$|\Sha|$,
the percentage increases from 3.6\% to 5.9\% to 8.0\%.
It seems unlikely that these percentages (for $p\neq 3$)
will climb all the way to~100\%, and without a better guess,
one could posit that they are tending toward the number suggested
by the Delaunay heuristic.
In Table~\ref{tbl:first}, the ``$r>0$'' column
counts percentages of curves for which the central $L$-value vanishes,
while the other four columns denote how often a given prime divides
the $|\Sha|$-value of a nonvanishing twist.

\begin{table}[hbtp]
\caption{Data for cubic twists}\label{tbl:first}
\begin{center}
\begin{tabular}{|c|c|c|c|c|c|}\hline
&$r>0$&$p=2$&$p=3$&$p=5$&$p=7$\\\hline
$m\le 10^5$&22.9&37.3&33.7&3.9&1.2\\
$m\le 10^6$&20.5&42.1&40.1&5.9&2.4\\
$m\le 10^7$&17.7&45.7&45.3&8.0&3.7\\
prime $m\le 10^7$&6.1&53.5&5.8&14.5&8.2\\
Delaunay& &58.3&36.1&20.7&14.5\\\hline
\end{tabular}
\end{center}
\end{table}

It was pointed out to us by M.~O.~Rubinstein that quadratic twist data for
$|\Sha|$ tend to the Delaunay number more readily upon including {\it all}
even rank twists, instead of just the ones of rank~0. Indeed, as we expect
that the high rank twists should form an asymptotically negligible set,
there is perhaps no reason not to include them in our data. Furthermore,
additionally restricting to prime twists also tends to speed convergence
toward the number given by Delaunay. Upon implementing these two ideas,
we get numbers of $56.3\%$ for 2-divisibility,
$19.7\%$ for 5-divisibility, and $13.8\%$ for 7-divisibility,
which are fairly close to the percentages predicted by Delaunay.
For 3-divisibility we have only $11.6\%$, as the existence of
\hbox{3-isogenies} for our curves appears to have a definite impact
(Rubinstein reports similar phenomena for quadratic twists).

One can do a similar experiment with quartic twists of~$X_0(32)$
or sextic twists of $X_0(27)$. We only looked at the former.
For the computation of the sign of the functional equation in these cases,
see~\cite{st}. Note that \cite{zk} look at the quadratic twists of $X_0(32)$
given by \hbox{$y^2=x^3-m^2x$} with $m\equiv 1$ (mod~16) for~$m\le 500000$,
and they find that the percentage of vanishing twists is dropping fairly
rapidly, it being 15.2\% for $m\le 50000$ and 10.6\% for $m\le 500000$.
For the quartic twists of $X_0(32)$ we are looking at $y^2=x^3+mx$
where 4 does not divide $m$ and $m$ is free of fourth powers.
Here we consider positive $m\le 8000000$,
of which 24.9\% of the even twists have vanishing central value.
This is less than the 27.4\% for $m\le 10^6$, and 29.8\% for $m\le 10^5$.
Similar percentages occur for the negative~$m$.

\section{Computational techniques}
The computations were carried out on a network of about 10 SPARC machines
(mostly SPARC-V) over a 6-month period at the beginning of~2001.
Our bound of $m\le 10^7$ was chosen as we were mainly interested in the
question of extra vanishing, and seemed sufficient to answer the question
posed by \cite{zk} on whether the rate remained constant.
With today's technology, extending the experiment to $m\le 10^8$
should be feasible, as should a similar experiment looking at
cubic twists with odd functional equation.

As stated in~\cite{zk}, the computation of the $a_m(n)$
takes time $O(\log n)$ if $n$ is prime and $O(1)$ time otherwise
(using the multiplicativity relations, viewing the values for the
primes dividing~$n$ as taking negligible time as they are already computed).
We computed the values of $a_1(p)$ for $p\le 10^9$ once-and-for-all
ahead of time, and then read these from disk as needed.
Additionally, tricks such as fast modular exponentiation were used
to speed up the computation of $m^{(p-1)/3}$ mod~$p$.
Similarly, the computing of $e^{-2\pi n/\sqrt N}$
was faciliated by the fact that the exponential function is
a homomorphism; for a given~$N$, we computed various powers
of $e^{-2\pi/\sqrt N}$ and then for each~$n$ multiplied these
together as needed to get the desired value.
For the computation of $L(E_m,1)$, and the question of how far
the infinite sum need be computed, we followed a method similar
to that of~\cite{zk}, calculating the $|\Sha|$-value
$S_m=\frac{T^2}{c\Omega}L(E_m,1)$
where $T$ is the size of the torsion group, $c$ is the global Tamagawa number,
and $\Omega$ is the real period (see pages 54--56 of \cite{zk} for these).
We then stop the calculation when $S_m$ is sufficiently close to an integer
(possibly zero). As a check, we expect all the $S_m$ values to be squares,
which indeed does turn out to be the case.

\section{Random matrix models}
In this section we make some comments about random matrix theory and
the expected number of even cubic twists of $X_0(27)$ which have vanishing
central $L$-value. We follow the ideas of \cite{ckrs} and~\cite{dfk}.
In our case of cubic twists, we expect, as do~\cite{ckrs}, to have
symmetry type~$O^+$, that is, orthogonal with positive determinant.
This is because the sign of our functional equation is always~$+1$.
Note that \cite{dfk} have unitary symmetry in their type of cubic twist,
due to the fact that the functional equation has an essentially arbitrary
complex number (related to a Gauss sum) appearing in~it.

We write $E=X_0(27)$ and $E_d$ for the $d$th cubic twist of~$E$.
As given in equations (20), (22), and (16) of~\cite{ckrs}, the assumption of
$O^+$ symmetry implies that $P_E(N,x)=c_E N^{3/8}/\sqrt x$
should approximate (for small~$x$) the probability density function
for values of $L(E_d,1)$, where $N\sim \log X$ and we integrate
$\int_0^X P_E(N,x)\, dx$ to get an expected probability that
$L(E_d,1)$ is less than~$X$. The idea is that we know that the actual
values of $L(E_d,1)$ are discretised (due to the Birch--Swinnerton-Dyer
formula), and thus we declare (in a somewhat arbitrary manner) sufficiently
small values of $L(E_d,1)$ to indicate that in fact we have $L(E_d,1)=0$.
We recall that BSD implies we have
$$\frac{L(E_d,1)}{\Omega_d}=\prod_{p|3d} c_p\cdot {\frac{|\Sha_d|}{|T_d|^2}}$$
where $\Omega_d$ is the real period of~$E_d$, the $c_p$ are Tamagawa numbers,
$\Sha_d$ is the Shafarevitch--Tate group, and $T_d$ is the torsion
group of~$E_d$. We are thus thinking of $|\Sha_d|$ (which is a square)
as our discretised variable, with everything else being computable.
When $d>2$ the torsion group is trivial.
For cubefree~$d$ we have that $\Omega_d=\Omega_1/d^{1/3}$,
except when $9|d$ in which case we have $\Omega_d=3\Omega_1/d^{1/3}$.
Note that in definition (8) of~\cite{ckrs}, quadratic twists that are
not relatively prime to the conductor are excluded; we will similarly
exclude twists that are divisible by~3, though one could deal with them
via making appropriate corrections. For the Tamagawa product we have that
$c_3=3$ when $d\equiv\pm 1$ (mod~9), $c_3=2$ when $d\equiv\pm 2$ (mod~9),
and $c_3=1$ otherwise, while $c_p=3$ for primes $p\equiv 1$ (mod~3)
and $c_p=1$ for primes $p\equiv 2$ (mod~3). Given this divergent behaviour
based upon prime divisibility, as in Conjecture~1 of \cite{ckrs} we decided
to restrict to prime twists, and additionally split the primes into
congruence classes modulo~9. Indeed, it is calculable that the sign
of the functional equation is odd when our prime twist $d$ is congruent
to $4,7,8$ (mod~9), and by 3-descent we can verify that the rank is zero
when $d$ is 2 or 5 (mod~9). Moreover, again by 3-descent,
we know that the rank is at most 2 (and the functional equation is even)
when $d$ is 1~mod~9.
Computing as with equation (23) in \cite{ckrs} we get the following:

\begin{question}
Let $V_T$ be the set of primes $d$ less than $T$
congruent to 1 modulo~9 with $L(E_d,1)=0$.
Is there some constant $c\neq 0$ such that
$$\sum_{d\in V_T} 1\sim cT^{5/6}(\log T)^{-5/8}\qquad?$$
\end{question}

Assuming an affirmative answer,
our data give a constant of approximately $c=1/6$.
The argument is similar for quartic twists of $X_0(32)$
or sextic twists of $X_0(27)$, and we can expect asymptotics
for prime twists of order $T^{7/8}(\log T)^{-5/8}$
and $T^{11/12}(\log T)^{-5/8}$, and upon restricting to various congruence
classes we should get appropriate constants in front of these.
Via techniques from prime number theory and considerations
from Tamagawa numbers, one should be able to argue as in \cite{ckrs}
to get an asymptotic for all cubefree twists.

Finally we derive a version of Conjecture~2 of \cite{ckrs} suitable
for cubic, quartic, and sextic twists. For cubic twists,
for a given prime \hbox{$p\equiv 1$ (mod~3)} there are
3 solutions to $a^2+3b^2=4p$ with~\hbox{$a\equiv 2$ (mod~3)},
which correspond to the three possibilities for
the Frobenius trace~$a_p$. The argument given from (27)-(31)
in \cite{ckrs} does not differ (see below),
and so we get the following:

\begin{question}
Let $p\ge 5$ be prime, and for $1\le q\le p-1$ let $F_p^q(T)$ be the set
of cubefree positive integers $d\equiv q$ (mod~$p$) that are less than~$T$
such that $x^3+y^3=d$ has even functional equation.
Letting $a_d(p)$ be the $p$th trace of Frobenius for $x^3+y^3=d$
(where $d$ need not be cubefree), do we have
$$\lim_{T\rightarrow\infty}
\biggl(\sum_{d\in F_p^Y(T)} 1\biggm/
\sum_{d\in F_p^Z(T)} 1\biggr)=
\sqrt{\frac{p+1-a_Y(p)}{p+1-a_Z(p)}}\qquad?$$
\end{question}

We can also make a similar calculation for quartic and sextic twists.
In Tables \ref{tab:275}-\ref{tab:2711}
below we list vanishing probabilities in support of
an affirmative answer to the above question;
the $c$-column represents which congruence class is used.
For $p=7$ the ratios should be $[\sqrt{3}:\sqrt{9}:\sqrt{12}]$,
and for $p=13$ they should be $[\sqrt{9}:\sqrt{12}:\sqrt{21}]$.

We also have some data (see Tables \ref{tab:325}-\ref{tab:3211})
for the vanishing frequencies for positive quartic twists of~$X_0(32)$.
For $p=5$ the ratios should be given by
$\bigl[\sqrt{2}:\sqrt{4}:\sqrt{8}:\sqrt{10}\bigr]$;
for $p=13$ they should be
$\bigl[\sqrt{8}:\sqrt{10}:\sqrt{18}:\sqrt{20}\bigr]$.

The heuristic for Conjecture~2 in~\cite{ckrs} is based upon supposed
cancellation from a quadratic character, whereas in our cubic twist case
the source of cancellation is perhaps not so transparent.
Therefore we go through the details.
We have that
$$\sum_{d\in F_p^q(T)} L(E_d,1/2)^k=
\sum_{d\in F_p^q(T)}\biggl(\sum_{n=1}^\infty \frac{a_d(n)}{n}\biggr)^k=
\sum_{d\in F_p^q(T)}\sum_{n=1}^\infty \frac{b_d(n)}{n},$$
where $b_m(n)=\sum_{n=n_1\cdots n_k} a_m(n_1)\cdots a_m(n_k)$
with the sum being over all ways of writing $n$
as a product of $k$ positive factors.
If we invert the order of summation in this last expression,
the sum over $d$ should typically have much cancellation
since the $b_d(n)$ are essentially randomly distributed.
This, however, is not the case for $n$ that are a power of~$p$,
as here the value of $a_d(p^r)$ is fixed since $d$ is fixed modulo~$p$.
Thus we should get a main contribution in the above by restricting to
values of $n$ that are powers of~$p$ (indeed, if we did this argument with
no congruence restriction we would expect $n=1$ to give the main term).
As in (31) of \cite{ckrs} we thus get that
\begin{align*}
\sum_{d\in F_p^q(T)} L(E_d,1/2)^k\sim
&\sum_{d\in F_p^q(T)}\sum_{p^r} \frac{b_d(p^r)}{p^r}=
\sum_{d\in F_p^q(T)}\biggl(\sum_{p^r} \frac{a_d(p^r)}{p^r}\biggr)^k=\\
&\qquad=\biggl(\frac{p}{p+1-a_d(p)}\biggr)^k\sum_{d\in F_p^q(T)} 1.\\
\end{align*}
We complete our heuristic by first noting that the
sets $F_p^q(T)$ have asymptotically equal sizes and then taking $k=-1/2$
as is suggested by the random matrix theory of~\cite{ckrs}.
Note that a similar heuristic can be given for moments of higher
derivatives, but the combinatorics become more difficult due to
the presence of logarithms. In this context, the data of
Elkies \cite{elkies-antsv} distinctly show a congruence-class
phenomenon for rank 3 quadratic twists of~$X_0(32)$.

\section{Acknowledgments}
The author was partially funded by an NSF VIGRE Postdoctoral Fellowship
at The Pennsylvania State University for part of the time this work
was done. He also thanks an anonymous referee for useful comments.

\def\foreskip{\vspace{-12pt}}
\def\aftskip{\vspace{-7pt}}
\begin{table}[hbtp]
\begin{minipage}{2in}
\caption{$p=5$, $X_0(27)$\label{tab:275}}
\begin{center}
\begin{tabular}{|c|c|c|c|}\hline
$c$&$\#r>0$&\#curves&\\\hline
%0&173842&805068&0.216\\
1&140463&838612&0.167\\
2&140549&838570&0.168\\
3&140613&838575&0.168\\
4&140750&838637&0.168\\\hline
\end{tabular}
\end{center}
\foreskip\caption{$p=7$, $X_0(27)$}\aftskip
\begin{center}
\begin{tabular}{|c|c|c|c|}\hline
$c$&$\#r>0$&\#curves&\\\hline
%0&146882&583813&0.252\\
1&109569&595982&0.184\\
2&125728&595952&0.211\\
3&59440&595912&0.100\\
4&58759&595903&0.099\\
5&125714&595963&0.211\\
6&110125&595937&0.185\\\hline
\end{tabular}
\end{center}
\foreskip\caption{$p=11$, $X_0(27)$\label{tab:2711}}\aftskip
\begin{center}
\begin{tabular}{|c|c|c|c|}\hline
$c$&$\#r>0$&\#curves&\\\hline
%0&85515&375291&0.228\\
1&64989&378410&0.172\\
2&65211&378408&0.172\\
3&65001&378430&0.172\\
4&65008&378444&0.172\\
5&64956&378423&0.172\\
6&65208&378426&0.172\\
7&65054&378411&0.172\\
8&64773&378422&0.171\\
9&65164&378396&0.172\\
10&65338&378401&0.173\\\hline
\end{tabular}
\end{center}
\foreskip\caption{$p=13$, $X_0(27)$\label{tab:2713}}\aftskip
\begin{center}
\begin{tabular}{|c|c|c|c|}\hline
$c$&$\#r>0$&\#curves&\\\hline
%0&81353&318223&0.256\\
1&44504&320075&0.139\\
2&52214&320099&0.163\\
3&51754&320124&0.162\\
4&67352&320151&0.210\\
5&43064&320116&0.135\\
6&68325&320090&0.213\\
7&68702&320124&0.215\\
8&43215&320104&0.135\\
9&67584&320107&0.211\\
10&51465&320072&0.161\\
11&51827&320135&0.162\\
12&44858&320042&0.140\\\hline
\end{tabular}
\end{center}
\end{minipage}
\begin{minipage}{2in}
\caption{$p=5$, $X_0(32)$\label{tab:325}}
%\vspace{0pt}
\begin{center}
\begin{tabular}{|c|c|c|c|}\hline
$c$&$\#r>0$&\#curves&\\\hline
%0&218095&743145&0.293\\
1&156097&749089&0.208\\
2&104136&749107&0.139\\
3&236861&749125&0.316\\
4&215944&749182&0.288\\\hline
\end{tabular}
\end{center}
\foreskip\caption{$p=7$, $X_0(32)$}\aftskip
\begin{center}
\begin{tabular}{|c|c|c|c|}\hline
$c$&$\#r>0$&\#curves&\\\hline
%0&160258&508591&0.315\\
1&128846&538523&0.239\\
2&128491&538505&0.239\\
3&128553&538517&0.239\\
4&128597&538505&0.239\\
5&128053&538495&0.238\\
6&128335&538512&0.238\\\hline
\end{tabular}
\end{center}
\foreskip\caption{$p=11$, $X_0(32)$\label{tab:3211}}\aftskip
\begin{center}
\begin{tabular}{|c|c|c|c|}\hline
$c$&$\#r>0$&\#curves&\\\hline
%0&104754&328844&0.319\\
1&82653&341092&0.242\\
2&82782&341070&0.243\\
3&82581&341069&0.242\\
4&82392&341072&0.242\\
5&82806&341113&0.243\\
6&82448&341061&0.242\\
7&82661&341108&0.242\\
8&82388&341045&0.242\\
9&82720&341091&0.243\\
10&82948&341083&0.243\\\hline
\end{tabular}
\end{center}
\foreskip\caption{$p=13$, $X_0(32)$\label{tab:3213}}\aftskip
\begin{center}
\begin{tabular}{|c|c|c|c|}\hline
$c$&$\#r>0$&\#curves&\\\hline
%0&86777&287544&0.302\\
1&85079&287669&0.296\\
2&60843&287670&0.212\\
3&85408&287673&0.297\\
4&53551&287693&0.186\\
5&60788&287689&0.211\\
6&60926&287684&0.212\\
7&81716&287656&0.284\\
8&81500&287704&0.283\\
9&85480&287661&0.297\\
10&53852&287654&0.187\\
11&81525&287683&0.283\\
12&53688&287668&0.187\\\hline
\end{tabular}
\end{center}
\end{minipage}
\end{table}

\end{document}